# Cluster soft sets and cluster soft topologies

Zanyar A. Ameen [1] and Samer Al Ghour [2*]
[1]Department of Mathematics, College of Science, University of Duhok,
Duhok 42001, Iraq; ORCID 0000-0003-0740-3331
[2]Department of Mathematics and Statistics, Jordan University of Science and Technology,
Irbid 22110, Jordan; algore@just.edu.jo; ORCID 0000-0002-7872-440X

The cluster soft point is an attempt to introduce a novel generalization of the soft closure point and the soft limit point. A cluster soft set is defined to be the system of all cluster soft points of a soft set. Then the fundamental properties of cluster soft sets are demonstrated. Moreover, the concept of a cluster soft topology on a universal set is introduced with regard to the cluster soft sets. The cluster soft topology is derived from a soft topology with an associated soft ideal, but it is finer than the original soft topology. On the other hand, if we start constructing the cluster soft topology from another cluster soft topology, we will end up with the first cluster soft topology we started with. The implication of cluster soft topologies is highlighted using some examples. Eventually, we represent the cluster soft closed sets in terms of several forms of soft sets.

## 1 Introduction

For researchers in various disciplines, including medical science, economics, science, and engineering, uncertainty and insufficient information present numerous challenges. Soft set, defined by Molodtsov (1999), is a significant direction that has given rise to several extensions for overcoming these challenges among the different approaches intended to address them (see Molodtsov (2004)). Other theories, such as rough set theory of Pawlak (1982) and fuzzy set theory of Zadeh (1965), can be considered as mathematical methods for dealing with uncertainty, but each has its own set of difficulties. Ali et al. (2009) developed soft set theory by defining some new operations. The structure of parameter sets, especially those linked to soft sets, provides a consistent framework for modeling uncertain data. This results in the fast development of soft set theory in a short period of time, as well as a large variety of soft set real-world applications (see Al-shami (2022); Dalkılıç¸ (2021, 2022b); Dalkılıç and Demirtaş (2022); Liu et al. (2021). For more details, we refer the readers to the survey given by Zhan and Alcantud (2019).

It is known that rough set theory uses equivalence classes, while fuzzy set theory depends on the grades of memberships. In soft set theory, membership is determined by appropriate parameters. Despite being very different from one another, all three theories address uncertainty. An effective outcome from the mixing of these theories is possible. Such theories nowadays are called hybrid methods or models. Dubois and Prade (1990) were the first to combine fuzzy and rough set theories. In this fashion many generalization of soft sets appeared in the literature. For example, fuzzy soft set and rough soft set theories were respectively established by Maji et al. (2001) and Ali (2011). For more known and recent types of hybrid models, we encourage readers to check (Santos-García and Alcantud (2023)). These models are not free from the real applications (see Dalkılıç (2022a); Maji et al. (2001); Santos-García and Alcantud (2023)).

Multiple researchers have applied soft set theory to various mathematical structures like soft group theory Aktaş and Çağman (2007), soft ring theory Acar et al. (2010), soft category theory Sardar and Gupta (2013), etc. Soft topology is one of the structures introduced by Shabir and Naz (2011) and Çağman et al. (2011) as a novel generalization of the classical topology. The work in the latter two manuscripts was crucial to the development of the field of soft topology. Many classical concepts in topology have been generalized and extended in soft set settings, for instance, soft separation axioms(Al-shami and El-Shafei (2020); Bayramov and Gunduz (2018)), soft separable spaces (Bayramov and Gunduz (2018)), soft connected spaces (Lin (2013)), soft compact spaces (Aygünoğlu and Aygün, soft paracompact spaces (Lin (2013)), soft extremally disconnected spaces (Asaad (2017)), and soft submaximal spaces (Al Ghour and Ameen (2022)). The maximal or minimal elements concerning certain soft topological properties in the lattice of all soft topologies have been studied in (Al Ghour and Ameen (2022a); Ameen and Al Ghour (2022a)).

Methods of generating soft topologies over a common universal set are another useful line of research. Terepeta (2019) gave two remarkable formulas for generating soft topologies from crisp topologies. Then, Al-shami and Kocinac (2019) showed that the soft topology generated via one of the formulas is equivalent to the enriched soft topology. Alcantud (2020) improved the formulas given by Terepeta in such a way that one can generate a soft topology from a system of crisp topologies. Ameen and Al Ghour (2022) introduced the so-called soft simple extension of a soft topology. The simple extended soft topology with respect to a soft topology and a soft set is generated by their (soft) union. Azzam et al. (2022) introduced methods of generating soft topologies by various soft set operators. Among them, the soft closure and soft derived set operators. The cluster soft set can be seen as a generalized version of the latter soft set operators. The newly produced topology, from cluster soft sets, is called the cluster soft topology. Cluster soft topological spaces are a mix of soft topological and soft algebraic structures. That is to say, the cluster soft topology is obtained from a soft topology along with a soft ideal. By virtue of some examples, it is shown that cluster soft topologies are the most natural class of soft topologies. Furthermore, the category of cluster soft topologies includes the earlier mentioned soft topologies.

The rest of the paper is organized as follows. In Section 2, we provide a summary of the literature on soft set theory and soft topology. In Section 3, we introduce the notion of cluster soft sets and their fundamental properties to establish a new soft topology. In Section 4, the concept of cluster soft topology is introduced. Moreover, it is determined how original soft topology and cluster soft topology relate to one another. In Section 5, we propose characterizations of cluster soft closed sets in terms of some different types of soft sets. Finally, we conclude our paper with a short summary and possible lines of future work.

## 2 Preliminaries

Let $X$ be a universal set, $\widetilde{\Omega}$ a set of parameters, and $\mathcal{P}(X)$ the set of all subsets of $X$. If $F\colon \Omega \to \mathcal{P}(X)$ is a set-valued mapping and $\Omega \subseteq \widetilde{\Omega}$, then the collection $(F, \Omega) = \{(\alpha, F(\alpha))\colon \alpha \in \Omega\}$ is said to be a soft set over $X$. By $S_{\Omega}(X)$ we mean the set of soft sets over $X$ parameterized by $\Omega$. The soft complement $(F,\Omega)^c$(Shabir and Naz (2011)) of a soft set $(F,\Omega)$ is a soft set $(F^c,\Omega)$, where $F^c\colon \Omega \to \mathcal{P}(X)$ is a mapping having the property that $F^c(\alpha) = X - F(\alpha)$ for all $\alpha \in \Omega$. A soft set $(F,\Omega) \in S_{\Omega}(X)$ is called null (Maji et al. (2003)), denoted by $\widetilde{\Phi}$, if $F(\alpha) = \emptyset$ for all $\alpha \in \Omega$ and called absolute (Maji et al. (2003)), denoted by $\tilde{X}$, if $F(\alpha) = X$ for all $\alpha \in \Omega$. Evidently, $\tilde{X}^c = \widetilde{\Phi}$ and $\widetilde{\Phi}^c = \tilde{X}$. A soft set $(F,\Omega)$ is called finite (Das and Samanta (2013)) if $F(\alpha)$ is finite for each $\alpha \in \Omega$. Otherwise, it is called infinite. A soft element (Nazmul and Samanta (2013)), denoted by $x_\alpha$, is a soft set $(F,\Omega)$ over $X$ whenever $F(\alpha) = \{x\}$ and $F(\lambda) = \emptyset$ for all $\lambda \in \Omega$ with $\lambda \neq \alpha$, where $\alpha \in \Omega$ and $x \in X$. The soft element is called a soft point in (Xie (2015)). We prefer to use the concept of soft point in the sequel. By a statement $x_\alpha \in (F,\Omega)$ we mean $x \in F(\alpha)$. By $P_{\Omega}(X)$ we denote the set of all soft points in $X$. A soft set $(A,\Omega_1)$ is a soft subset of $(B,\Omega_2)$ (see Molodtsov (1999)), write $(A,\Omega_1) \widetilde{\subseteq} (B,\Omega_2)$, if $\Omega_1 \subseteq \Omega_2 \subseteq \widetilde{\Omega}$ and $A(\alpha) \subseteq B(\alpha)$ for all $\alpha \in \Omega_1$, and $(A,\Omega_1) = (B,\Omega_2)$ if $(A,\Omega_1) \widetilde{\subseteq} (B,\Omega_2)$ and $(B,\Omega_2) \widetilde{\subseteq} (A,\Omega_1)$. The soft union of soft sets $(A,\Omega),(B,\Omega)$ is represented by $(F,\Omega) = (A,\Omega) \mathbin{\widetilde{\cup}} (B,\Omega)$, where $F(\alpha) = A(\alpha) \cup B(\alpha)$ for all $\alpha \in \Omega$. The soft intersection of soft sets $(A,\Omega),(B,\Omega)$ is given by $(F,\Omega) = (A,\Omega) \mathbin{\widetilde{\cap}} (B,\Omega)$, where $F(\alpha) = A(\alpha) \cap B(\alpha)$ for all $\alpha \in \Omega$. The soft set difference $(A,\Omega) - (B,\Omega)$ is the soft set $(F,\Omega) = (A,\Omega) - (B,\Omega)$, where $F(\alpha) = A(\alpha) - B(\alpha)$ for $\alpha \in \Omega$ (see Terepeta (2019)).

The definitions of soft union and soft intersection of two soft sets with respect to arbitrary subsets of $\Omega$ were given by Maji et al. (2003). But it turns out that these definitions are misleading and ambiguous as reported by Ali et al. (2009) and Terepeta (2019). We follow that the definitions of soft union and soft intersection of soft sets given by Terepeta, which coincide with that adopted by Ali at al. (2009).

**Definition 2.1** (Çağman et al. (2011); Shabir and Naz (2011)) *A collection $\tilde{\mathcal{T}}$ of $S_{\Omega}(X)$ is said to be a soft topology on $X$ if it satisfies the following axioms:*

1. $\widetilde{\Phi}, \tilde{X} \in \tilde{\mathcal{T}}$.
2. If $(F,\Omega),(G,\Omega) \in \tilde{\mathcal{T}}$, then $(F,\Omega) \mathbin{\widetilde{\cap}} (G,\Omega) \in \tilde{\mathcal{T}}$.
3. If $\{(F_i,\Omega)\colon i \in I\} \widetilde{\subseteq} \tilde{\mathcal{T}}$, then $\widetilde{\bigcup}_{i\in I} (F_i,\Omega) \in \tilde{\mathcal{T}}$.

We call the triple $(X,\tilde{\mathcal{T}},\Omega)$ a soft topological space on $X$. We call the elements of $\tilde{\mathcal{T}}$ soft open sets. We call the complement of every soft open or elements of $\tilde{\mathcal{T}}^c$ soft closed sets. By $T_{\Omega}(X)$ we mean the lattice of all soft topologies on $X$ (see Al Ghour and Ameen (2022a)).

**Definition 2.2** (Nazmul and Samanta (2013)) *Let $(N,\Omega) \in S_{\Omega}(X)$ and $\tilde{\mathcal{T}} \in T_{\Omega}(X)$. Then $(N,\Omega)$ is called a soft neighborhood of $x_\alpha \in P_{\Omega}(X)$ if there exists $(U,\Omega) \in \tilde{\mathcal{T}}(x_\alpha)$ such that $x_\alpha \in (U,\Omega) \widetilde{\subseteq} (N,\Omega)$, where $\tilde{\mathcal{T}}(x_\alpha)$ is the family of all elements of $\tilde{\mathcal{T}}$ that contain $x_\alpha$.*

**Definition 2.3** (Çağman et al. (2011)) *Given a soft topology $\tilde{\mathcal{T}}$. A (countable) soft base for $\tilde{\mathcal{T}}$ is a (countable) subcollection $\mathcal{B} \subseteq \tilde{\mathcal{T}}$ such that elements of $\tilde{\mathcal{T}}$ are unions of elements of $\mathcal{B}$.*

**Definition 2.4** (Ameen and Al Ghour (2022a)) *Let $\mathcal{F} \widetilde{\subseteq} S_{\Omega}(X)$. The intersection of all soft topologies on $X$ containing $\mathcal{F}$ is called a soft topology generated by $\mathcal{F}$ and is referred to $\tilde{\mathcal{T}}[\mathcal{F}]$.*

**Lemma 2.5** (Shabir and Naz (2011)) *Let $(X,\tilde{\mathcal{T}},\Omega)$ be a soft topological space, then for each $\alpha \in \Omega$, the collection $\tilde{\mathcal{T}}(\alpha) = \{F(\alpha)\colon (F,\Omega) \in \tilde{\mathcal{T}}\}$ is a (crisp) topology on $X$.*

**Definition 2.6** (Shabir and Naz (2011)) *Let* $(B,\Omega) \in S_{\Omega}(X)$ *and* $\tilde{\mathcal{T}} \in T_{\Omega}(X)$.

1. The soft closure of $(B,\Omega)$ is $cl(B,\Omega) := \tilde{\cap} \{(F,\Omega): (B,\Omega) \tilde{\subseteq} (F,\Omega), (F,\Omega) \in \tilde{\mathcal{T}}^c\}$.
2. The soft interior of $(B,\Omega)$ is $int(B,\Omega) := \tilde{\cup} \{(F,\Omega): (F,\Omega) \tilde{\subseteq} (B,\Omega), (F,\Omega) \in \tilde{\mathcal{T}}\}$.

**Definition 2.7** (Çağman et al. (2011)) *Let* $(B,\Omega) \in S_{\Omega}(X)$ *and* $\tilde{\mathcal{T}} \in T_{\Omega}(X)$. *A soft point* $x_{\alpha} \in P_{\Omega}(X)$ *is called a limit soft point of* $(B,\Omega)$ *if* $(G,\Omega) \tilde{\cap} (B,\Omega) - \{x_{\alpha}\} \neq \tilde{\Phi}$ *for all* $(G,\Omega) \in \tilde{\mathcal{T}}(x_{\alpha})$. *The set of all limit soft points is symbolized by* $\mathcal{D}(B,\Omega)$.

Then $cl(F,\Omega) = (F,\Omega) \tilde{\cup} \mathcal{D}(F,\Omega)$ (see Theorem 5 in Çağman et al. (2011)).

**Definition 2.8** (Kandil et al. (2014)) *A non-null class* $\tilde{I} \tilde{\subseteq} S_{\Omega}(X)$ *is termed a soft ideal on* $X$ *if* $\tilde{I}$ *satisfies the following conditions:*

1. If $(R,\Omega),(S,\Omega) \in \tilde{I}$, then $(R,\Omega) \tilde{\cup} (S,\Omega) \in \tilde{I}$; and
2. If $(R,\Omega) \in \tilde{I}$ and $(S,\Omega) \tilde{\subseteq} (R,\Omega)$, then $(S,\Omega) \in \tilde{I}$.

$\tilde{I}$ is called a soft $\sigma$-ideal if (1) holds for countably many soft sets. We denote the family of soft ideals on $X$ by $I_{\Omega}(X)$.

Applying the definition, for any $\tilde{I},\tilde{J} \in I_{\Omega}(X)$, one can directly show that $\tilde{I} \tilde{\cap} \tilde{J}$ and $\tilde{I} \tilde{\Cup} \tilde{J}$ are also soft ideals, where $\tilde{I} \tilde{\Cup} \tilde{J} := \{(A,\Omega) \tilde{\cup} (B,\Omega): (A,\Omega) \in \tilde{I}, (B,\Omega) \in \tilde{J}\}$.

**Lemma 2.9** (Matejdes (2016); Matejdes (2021)) *Let* $Gr(F) = \{(\alpha,x) \in \Omega \times X: x \in F(\alpha)\}$, *which is the graph of the set-valued function* $F: \Omega \to \mathcal{P}(X)$. *Then*

1. If $(\Omega \times X, \mathcal{T})$ is a topological space, then $(X, \tilde{\mathcal{T}}, \Omega)$ is a soft topological space, where $\tilde{\mathcal{T}} = \{(F,\Omega): Gr(F) \in \mathcal{T}\}$.
2. If $(X, \tilde{\mathcal{T}}, \Omega)$ is a soft topological space, then $(\Omega \times X, \mathcal{T})$ is a topological space, where $\mathcal{T} = \{Gr(F): (F,\Omega) \in \tilde{\mathcal{T}}\}$.

**Lemma 2.10** (Kandil et al. (2014)) *Let* $\tilde{I}$ *be a soft ideal on* $X$. *Then for each* $\alpha \in \Omega$, *the collection* $\tilde{I}(\alpha) = \{A(\alpha): (A,\Omega) \in \tilde{I}\}$ *is a crisp ideal on* $X$.

# 3 Cluster soft sets

**Definition 3.1** *Let* $(R,\Omega) \in S_{\Omega}(X)$, $\tilde{\mathcal{T}} \in T_{\Omega}(X)$, *and* $\tilde{I} \in I_{\Omega}(X)$. *A soft point* $x_{\alpha} \in P_{\Omega}(X)$ *is a cluster soft point of* $(R,\Omega)$ *if* $(R,\Omega) \tilde{\cap} (U,\Omega) \notin \tilde{I}$ *for each* $(U,\Omega) \in \tilde{\mathcal{T}}(x_{\alpha})$. *We shall not make difference between the terminologies "cluster soft point"and "soft cluster point". The set of all the cluster soft points of* $(R,\Omega)$ *is called the cluster soft set of* $(R,\Omega)$ *and is denoted by* $\mathfrak{c}_{(\tilde{\mathcal{T}},\tilde{I})}(R,\Omega)$ *or shortly* $\mathfrak{c}(R,\Omega)$.

**Remark 3.2** *Given* $(R,\Omega) \in S_{\Omega}(X)$, $\tilde{\mathcal{T}} \in T_{\Omega}(X)$, *and* $\tilde{I} \in I_{\Omega}(X)$. *We shall remark that*

1. If $\tilde{I} = \{\tilde{\Phi}\}$, then $\mathfrak{c}(R,\Omega) = cl(R,\Omega)$. That is, the soft cluster points are identical to the soft closure points of $(R,\Omega)$.

2. If $\tilde{I} = \{(F,\Omega): (F,\Omega) \in S_\Omega(X), (F,\Omega) \text{ is finite}\}$, then $\mathfrak{c}(R,\Omega) = \mathcal{D}(R,\Omega)$. That is, the soft cluster points are identical to the soft limit points of $(R,\Omega)$.

Now, we present some properties of soft cluster sets.

**Proposition 3.3** *Let $(R,\Omega),(S,\Omega) \in S_\Omega(X)$, $\tilde{\mathcal{T}} \in S_\Omega(X)$, and $\tilde{I} \in I_\Omega(X)$. The following properties hold:*

1. If $(R,\Omega) \in \tilde{I}$, then $\mathfrak{c}(R,\Omega) = \tilde{\Phi}$.
2. If $(R,\Omega) \tilde{\subseteq} (S,\Omega)$, then $\mathfrak{c}(R,\Omega) \tilde{\subseteq} \mathfrak{c}(S,\Omega)$.
3. $\mathfrak{c}((R,\Omega) \tilde{\cap} (S,\Omega)) \tilde{\subseteq} \mathfrak{c}(R,\Omega) \tilde{\cap} \mathfrak{c}(S,\Omega)$.
4. $\mathfrak{c}((R,\Omega) \tilde{\cup} (S,\Omega)) = \mathfrak{c}(R,\Omega) \tilde{\cup} \mathfrak{c}(S,\Omega)$.
5. $\mathfrak{c}(R,\Omega) - \mathfrak{c}(S,\Omega) \tilde{\subseteq} \mathfrak{c}((R,\Omega) - (S,\Omega))$.

***Proof***.

1. Let $(R,\Omega) \in \tilde{I}$ and $x_\alpha \in P_\Omega(X)$. Then for each $(U,\Omega) \in \tilde{\mathcal{T}}(x_\alpha)$, $(R,\Omega) \tilde{\cap} (U,\Omega) \tilde{\subseteq} (R,\Omega)$ implies $(R,\Omega) \tilde{\cap} (U,\Omega) \in \tilde{I}$. Therefore, $x_\alpha \notin \mathfrak{c}(R,\Omega)$ and so $\mathfrak{c}(R,\Omega) = \tilde{\Phi}$.

2. Let $x_\alpha \in P_\Omega(X)$. If $x_\alpha \notin \mathfrak{c}(S,\Omega)$, then there exists $(U,\Omega) \in \tilde{\mathcal{T}}(x_\alpha)$ such that $(S,\Omega) \tilde{\cap} (U,\Omega) \in \tilde{I}$. Since $(R,\Omega) \tilde{\subseteq} (S,\Omega)$, then $(R,\Omega) \tilde{\cap} (U,\Omega) \tilde{\subseteq} (S,\Omega) \tilde{\cap} (U,\Omega)$ and hence $(R,\Omega) \tilde{\cap} (U,\Omega) \in \tilde{I}$. Thus, $x_\alpha \notin \mathfrak{c}(R,\Omega)$.

3. Since $(R,\Omega) \tilde{\cap} (S,\Omega) \tilde{\subseteq} (R,\Omega)$ and $(R,\Omega) \tilde{\cap} (S,\Omega) \tilde{\subseteq} (S,\Omega)$. By (2), $\mathfrak{c}((R,\Omega) \tilde{\cap} (S,\Omega)) \tilde{\subseteq} \mathfrak{c}(R,\Omega)$ and $\mathfrak{c}((R,\Omega) \tilde{\cap} (S,\Omega)) \tilde{\subseteq} \mathfrak{c}(S,\Omega)$. Therefore, $\mathfrak{c}((R,\Omega) \tilde{\cap} (S,\Omega)) \tilde{\subseteq} \mathfrak{c}(R,\Omega) \tilde{\cap} \mathfrak{c}(S,\Omega)$.

4. By a similar technique of (3), one can obtain $\mathfrak{c}(R,\Omega) \tilde{\cup} \mathfrak{c}(S,\Omega) \tilde{\subseteq} \mathfrak{c}\big((R,\Omega) \tilde{\cup} (S,\Omega)\big)$. On the other hand, let $x_\alpha \in P_\Omega(X)$ and $x_\alpha \notin \mathfrak{c}(R,\Omega) \tilde{\cup} \mathfrak{c}(S,\Omega)$. Then there exist $(U,\Omega),(V,\Omega) \in \tilde{\mathcal{T}}(x_\alpha)$ such that $(R,\Omega) \tilde{\cap} (U,\Omega), (S,\Omega) \tilde{\cap} (V,\Omega) \in \tilde{I}$. So $(U,\Omega) \tilde{\cap} (V,\Omega) \in \tilde{\mathcal{T}}(x_\alpha)$ and $[(R,\Omega) \tilde{\cap} (U,\Omega)] \tilde{\cup} [(S,\Omega) \tilde{\cap} (V,\Omega)] \in \tilde{I}$. But $[(R,\Omega) \tilde{\cup} (S,\Omega)] \tilde{\cap} [(U,\Omega) \tilde{\cap} (V,\Omega)] \tilde{\subseteq} [(R,\Omega) \tilde{\cap} (U,\Omega)] \tilde{\cup} [(S,\Omega) \tilde{\cap} (V,\Omega)]$ implies $[(R,\Omega) \tilde{\cup} (S,\Omega)] \tilde{\cap} [(U,\Omega) \tilde{\cap} (V,\Omega)] \in \tilde{I}$. Hence, $x_\alpha \in \mathfrak{c}((R,\Omega) \tilde{\cup} (S,\Omega))$. The proof is finished.

5. Since $(R,\Omega) = ((R,\Omega) - (S,\Omega)) \tilde{\cup} ((R,\Omega) \tilde{\cap} (S,\Omega))$. By (4), $\mathfrak{c}(R,\Omega) = \mathfrak{c}((R,\Omega) - (S,\Omega)) \tilde{\cup} \mathfrak{c}((R,\Omega) \tilde{\cap} (S,\Omega))$. By (2), $\mathfrak{c}(R,\Omega) \tilde{\subseteq} \mathfrak{c}((R,\Omega) - (S,\Omega)) \tilde{\cup} \mathfrak{c}(S,\Omega)$. Therefore, $\mathfrak{c}(R,\Omega) - \mathfrak{c}(S,\Omega) \tilde{\subseteq} \mathfrak{c}((R,\Omega) - (S,\Omega))$.

For any index set $J$ and any family $\mathfrak{F}$ of finite subsets of $J$, we have the following statements:

**Proposition 3.4** *Let $(R_j,\Omega) \in S_\Omega(X)$ for $j \in J$, $\tilde{\mathcal{T}} \in T_\Omega(X)$, and $\tilde{I} \in I_\Omega(X)$. The following properties hold:*

1. $\mathfrak{c}(\tilde{\bigcup}_{j\in N}(R_j,\Omega)) = \tilde{\bigcup}_{j\in N}\mathfrak{c}(R_j,\Omega)$, where $N \in \mathfrak{F}$.
2. $\tilde{\bigcup}_{j\in J}\mathfrak{c}(R_j,\Omega) \tilde{\subseteq} \mathfrak{c}(\tilde{\bigcup}_{j\in J}(R_j,\Omega))$.

3. $\mathfrak{c}(\widetilde{\bigcap}_{j\in J}(R_j,\Omega)) \widetilde{\subseteq} \widetilde{\bigcap}_{j\in J}\mathfrak{c}(R_j,\Omega)$.

4. $\mathfrak{c}(\widetilde{\bigcup}_{j\in J}(R_j,\Omega)) = \widetilde{\bigcup}_{j\in J}\mathfrak{c}(R_j,\Omega)\ \widetilde{\cup}\ [\widetilde{\bigcap}_{N\in\mathfrak{F}}\mathfrak{c}(\widetilde{\bigcup}_{j\in J-N}(R_j,\Omega))]$.

***Proof*.**

1. Let $N$ be a finite subset of $J$. Since $(R_j,\Omega) \widetilde{\subseteq} \widetilde{\bigcup}_{j\in N}(R_j,\Omega)$ for each $j$, by (2) in Proposition 3.3, $\mathfrak{c}(R_j,\Omega) \widetilde{\subseteq} \mathfrak{c}(\widetilde{\bigcup}_{j\in N}(R_j,\Omega))$ and hence $\widetilde{\bigcup}_{j\in N}\mathfrak{c}(R_j,\Omega) \widetilde{\subseteq} \mathfrak{c}(\widetilde{\bigcup}_{j\in N}(R_j,\Omega))$. On the other hand, if there exists $x_\alpha \in P_\Omega(X)$ such that $x_\alpha \notin \widetilde{\bigcup}_{j\in N}\mathfrak{c}(R_j,\Omega)$, then, for each $j \in N$, there exists $(U_j,\Omega) \in \tilde{\mathcal{T}}(x_\alpha)$ such that $(R_j,\Omega)\ \widetilde{\cap}\ (U_j,\Omega) \in \tilde{I}$. Therefore, $\widetilde{\bigcup}_{j\in N}(U_j,\Omega) \in \tilde{\mathcal{T}}(x_\alpha)$ and $\widetilde{\bigcup}_{j\in N}[(R_j,\Omega)\ \widetilde{\cap}\ (U_j,\Omega)] \in \tilde{I}$. But

$$[\widetilde{\bigcup}_{j\in N}(R_j,\Omega)]\ \widetilde{\cap}\ [\widetilde{\bigcap}_{j\in N}(U_j,\Omega)] \widetilde{\subseteq} \widetilde{\bigcup}_{j\in N}[(R_j,\Omega)\ \widetilde{\cap}\ (U_j,\Omega)].$$

This implies that $[\widetilde{\bigcup}_{j\in N}(R_j,\Omega)]\ \widetilde{\cap}\ [\widetilde{\bigcup}_{j\in N}(U_j,\Omega)] \in \tilde{I}$ and so, $x_\alpha \notin \mathfrak{c}(\widetilde{\bigcup}_{j\in N}(R_j,\Omega))$. Thus, $\mathfrak{c}(\widetilde{\bigcup}_{j\in N}(R_j,\Omega)) \widetilde{\subseteq} \widetilde{\bigcup}_{j\in N}\mathfrak{c}(R_j,\Omega)$. Both of the inclusions prove (1).

2. Since $(R_j,\Omega) \widetilde{\subseteq} \widetilde{\bigcup}_{j\in J}(R_j,\Omega)$ for each $j \in J$, by (2) in Proposition 3.3, $\mathfrak{c}(R_j,\Omega) \widetilde{\subseteq} \mathfrak{c}(\widetilde{\bigcup}_{j\in J}(R_j,\Omega))$ and so $\widetilde{\bigcup}_{j\in J}\mathfrak{c}(R_j,\Omega) \widetilde{\subseteq} \mathfrak{c}(\widetilde{\bigcup}_{j\in J}(R_j,\Omega))$.

3. Since $\widetilde{\bigcap}_{j\in J}(R_j,\Omega) \widetilde{\subseteq} (R_j,\Omega)$ for each $j \in J$, by (2) in Proposition 3.3, $\mathfrak{c}\left(\widetilde{\bigcap}_{j\in J}(R_j,\Omega)\right) \widetilde{\subseteq} \mathfrak{c}(R_j,\Omega)$, for each $j$, and thus $\mathfrak{c}(\widetilde{\bigcap}_{j\in J}(R_j,\Omega)) \widetilde{\subseteq} \widetilde{\bigcap}_{j\in J}\mathfrak{c}(R_j,\Omega)$.

4. By (2), the inclusion $\widetilde{\bigcup}_{j\in J}\mathfrak{c}(R_j,\Omega)\ \widetilde{\cup}\ [\widetilde{\bigcap}_{N\in\mathfrak{F}}\mathfrak{c}(\widetilde{\bigcup}_{j\in J-N}(R_j,\Omega))] \widetilde{\subseteq} \mathfrak{c}(\widetilde{\bigcup}_{j\in J}(R_j,\Omega))$ can be followed. To prove the other direction, we choose any $N \in \mathfrak{F}$. From (1), we can have

$$\mathfrak{c}(\widetilde{\bigcup}_{j\in J}(R_j,\Omega)) = \widetilde{\bigcup}_{j\in N}\mathfrak{c}(R_j,\Omega)\ \ \widetilde{\cup}\ \ \mathfrak{c}(\widetilde{\bigcup}_{j\in J-N}(R_j,\Omega)).$$

Therefore,

$$\mathfrak{c}(\widetilde{\bigcup}_{j\in J}(R_j,\Omega)) \widetilde{\subseteq} \widetilde{\bigcup}_{j\in J}\mathfrak{c}(R_j,\Omega)\ \ \widetilde{\cup}\ \ \mathfrak{c}(\widetilde{\bigcup}_{j\in J-N}(R_j,\Omega)).$$

Since $N$ was chosen arbitrarily, so

$$\mathfrak{c}(\widetilde{\bigcup}_{j\in J}(R_j,\Omega)) \widetilde{\subseteq} \widetilde{\bigcup}_{j\in J}\mathfrak{c}(R_j,\Omega)\ \ \widetilde{\cup}\ \ [\widetilde{\bigcap}_{N\in\mathfrak{F}}\mathfrak{c}(\widetilde{\bigcup}_{j\in J-N}(R_j,\Omega))].$$

Hence the proof.

**Lemma 3.5** *Let $(R,\Omega) \in S_\Omega(X)$, $\tilde{\mathcal{T}} \in T_\Omega(X)$, and $\tilde{I} \in I_\Omega(X)$. Then*

1. $\mathfrak{c}(R,\Omega) \widetilde{\subseteq} cl(R,\Omega)$.

2. $\mathfrak{c}(R,\Omega) \in \tilde{\mathcal{T}}^c$.

3. $\mathfrak{c}[\mathfrak{c}(R,\Omega)] \widetilde{\subseteq} \mathfrak{c}(R,\Omega)$.

***Proof*.**

1. This follows as $\widetilde{\Phi} \in \tilde{I}$.

2. To prove (2) we relate each soft point $x_\alpha \in \tilde{X} - \mathfrak{c}(R,\Omega)$ to a soft set $(U_{x_\alpha},\Omega) \in \tilde{\mathcal{T}}(x_\alpha)$ such that $(R,\Omega) \tilde{\cap} (U_{x_\alpha},\Omega) \in \tilde{I}$. So for any $y_\lambda \in (U_{x_\alpha},\Omega)$, we have $y_\lambda \in \tilde{X} - \mathfrak{c}(R,\Omega)$. Then

$$\tilde{X} - \mathfrak{c}(R,\Omega) = \tilde{\bigcup}_{x_\alpha \in \tilde{X} - \mathfrak{c}(R,\Omega)} (U_{x_\alpha},\Omega).$$

Thus, $\tilde{X} - \mathfrak{c}(R,\Omega) \in \tilde{\mathcal{T}}$ as it is a union of soft open sets. Consequently, $\mathfrak{c}(R,\Omega) \in \tilde{\mathcal{T}}^c$.

3. By (1) and Proposition 3.3 (2), we have $\mathfrak{c}[\mathfrak{c}(R,\Omega)] \tilde{\subseteq} cl[\mathfrak{c}(R,\Omega)]$. But $cl[\mathfrak{c}(R,\Omega)] = \mathfrak{c}(R,\Omega)$ from (2). Hence, $\mathfrak{c}[\mathfrak{c}(R,\Omega)] \tilde{\subseteq} \mathfrak{c}(R,\Omega)$.

According to the above results, a soft cluster set is a generalization of a soft closure set and $(R,\Omega) \tilde{\subseteq} cl(R,\Omega)$ for any $(R,\Omega) \in S_\Omega(X)$. However, we cannot have $(R,\Omega) \tilde{\subseteq} \mathfrak{c}(R,\Omega)$ in general.

**Example 3.6** *For a set of parameters $\Omega = \{\alpha,\lambda\}$, let $\tilde{I} = \{(\alpha, A(\alpha)), (\lambda, A(\lambda)): A(\alpha) \subseteq \mathbb{Q}, A(\lambda)$ (finite) $\subseteq \mathbb{R}\}$ and let $\tilde{\mathcal{T}}$ be the soft topology on the set of real numbers $\mathbb{R}$ generated by*

$$\{\{(\alpha, B(\alpha)), (\lambda, B(\lambda))\}: B(\alpha) = (a,b), B(\lambda) = (c,d]; a,b,c,d \in \mathbb{R}; a < b, c < d\}.$$

Take $(R,\Omega) = \{(\alpha,\{5\}), (\lambda,\{\pi\})\}$. Then $\mathfrak{c}(R,\Omega) = \tilde{\Phi}$ and hence $(R,\Omega) \tilde{\nsubseteq} \mathfrak{c}(R,\Omega)$.

We end this part by commenting that the cluster soft set of $(R,\Omega)$ is called a soft local function $(^*)$ of $(R,\Omega)$ in Kandil et al. (2014). The reader may find some results in this direction in Kandil et al. (2014), but our theory is totally different. Furthermore, the conclusion (8) in Theorem 3.2 in Kandil et al. (2014) is false.

**Example 3.7** *If $\tilde{\mathcal{T}}$ is the soft topology on the set of real numbers $\mathbb{R}$ generated by*

$$\{((a,b),\Omega); a,b \in \mathbb{R}; a < b\},$$

where $\Omega$ is any set of parameters and $\tilde{I}$ is a soft ideal of finite soft subsets of $\tilde{\mathbb{R}}$. Take $(A_n,\Omega) = (\{1/n\},\Omega)$ for $n \in \mathbb{N}$. Then $\mathfrak{c}(A_n,\Omega) = \tilde{\Phi}$ for each $n$ and so $\tilde{\bigcup}_n \mathfrak{c}(A_n,\Omega) = \tilde{\Phi}$, while $\mathfrak{c}(\tilde{\bigcup}_n (A_n,\Omega)) = (\{0\},\Omega)$. Thus, $\mathfrak{c}(\tilde{\bigcup}_n (A_n,\Omega)) \neq \tilde{\bigcup}_n \mathfrak{c}(A_n,\Omega)$.

# 4 Cluster soft topologies

In order to introduce the cluster soft topology on $X$, we need to define the following concept:

**Definition 4.1** *Let $(R,\Omega) \in S_\Omega(X)$, $\tilde{\mathcal{T}} \in T_\Omega(X)$, and $\tilde{I} \in I_\Omega(X)$. Then $(R,\Omega)$ is said to be a cluster soft closed set (shortly, soft $\mathfrak{c}$-closed set) if $\mathfrak{c}(R,\Omega) \tilde{\subseteq} (R,\Omega)$.*

**Lemma 4.2** *Let $\tilde{\mathcal{T}} \in T_\Omega(X)$ and $\tilde{I} \in I_\Omega(X)$. The following statements are valid:*

1. $\tilde{\Phi}, \tilde{X}$ are soft $\mathfrak{c}$-closed.
2. Each element of $\tilde{I}$ is soft $\mathfrak{c}$-closed.
3. Each soft closed set is soft $\mathfrak{c}$-closed.
4. Any intersection of soft $\mathfrak{c}$-closed sets is soft $\mathfrak{c}$-closed.
5. A finite union of soft $\mathfrak{c}$-closed sets is soft $\mathfrak{c}$-closed.

***Proof***.

1. $\mathfrak{c}(\widetilde{\Phi}) = \widetilde{\Phi}$ by Proposition 3.3 (1) as $\widetilde{\Phi} \in \tilde{I}$ and $\mathfrak{c}(\tilde{X}) \widetilde{\subseteq} \tilde{X}$ always. Therefore, $\widetilde{\Phi}, \tilde{X}$ are soft $\mathfrak{c}$-closed.

2. Given $(R, \Omega) \in S_\Omega(X)$. If $(R, \Omega) \in \tilde{I}$, by Proposition 3.3 (1), $\mathfrak{c}(R, \Omega) = \widetilde{\Phi} \ \widetilde{\subseteq} (R, \Omega)$. Thus, $(R, \Omega)$ is soft $\mathfrak{c}$-closed.

3. Given $x_\alpha \in P_\Omega(X)$. Assume $(R, \Omega) \in \tilde{\mathcal{T}}^c$ and $x_\alpha \in \mathfrak{c}(R, \Omega)$. Then for each $(U, \Omega) \in \tilde{\mathcal{T}}(x_\alpha)$, $(R, \Omega) \widetilde{\cap} (U, \Omega) \neq \widetilde{\Phi}$. Since $(R, \Omega)^c \in \tilde{\mathcal{T}}$ and $(R, \Omega) \widetilde{\cap} (R, \Omega)^c = \widetilde{\Phi}$, we shall obtain that $x_\alpha \notin (R, \Omega)^c$. This implies that $x_\alpha \in (R, \Omega)$. Hence, $\mathfrak{c}(R, \Omega) \widetilde{\subseteq} (R, \Omega)$. This shows that $(R, \Omega)$ is soft $\mathfrak{c}$-closed.

4. Let $\{(R_j, \Omega): j \in J\}$ be a family of soft $\mathfrak{c}$-closed sets over $X$. For each $j$, we then have $\mathfrak{c}(R_j, \Omega) \widetilde{\subseteq} (R_j, \Omega)$. By Proposition 3.4 (3) and the latter statement, we have

$$\mathfrak{c}(\widetilde{\bigcap}_{j \in J} (R_j, \Omega)) \widetilde{\subseteq} \widetilde{\bigcap}_{j \in J} \mathfrak{c}(R_j, \Omega) \widetilde{\subseteq} \widetilde{\bigcap}_{j \in J} (R_j, \Omega).$$

Thus, $\widetilde{\bigcap}_{j \in J} (R_j, \Omega)$ is soft $\mathfrak{c}$-closed.

5. Let $(R_j, \Omega)$ be a soft $\mathfrak{c}$-closed sets over $X$, for $j = 1,2, \ldots, n$. By Definition 4.1, for each $j$, we have $\mathfrak{c}(R_j, \Omega) \widetilde{\subseteq} (R_j, \Omega)$ and therefore $\widetilde{\bigcup}_{j=1}^n \mathfrak{c}(R_j, \Omega) \widetilde{\subseteq} \widetilde{\bigcup}_{j=1}^n (R_j, \Omega)$. By Proposition 3.4 (1), we have

$$\mathfrak{c}(\widetilde{\bigcup}_{j=1}^{n} (R_j, \Omega)) = \widetilde{\bigcup}_{j=1}^{n} \mathfrak{c}(R_j, \Omega) \widetilde{\subseteq} \widetilde{\bigcup}_{j=1}^{n} (R_j, \Omega).$$

Hence, $\widetilde{\bigcup}_{j=1}^n (R_j, \Omega)$ is soft $\mathfrak{c}$-closed.

It is essential to note that the reverse of (3) is not true in general.

**Example 4.3** *Let $(X, \tilde{\mathcal{T}}, \Omega)$ be the indiscrete soft topological space and let $\tilde{I} = S_\Omega(X)$ be the soft ideal on $X$. Any proper soft subset of $\tilde{X}$ is soft $\mathfrak{c}$-closed but not soft closed.*

**Remark 4.4** *Given a soft topological space $(X, \tilde{\mathcal{T}}, \Omega)$ and a soft ideal $\tilde{I}$ on $X$. A soft subset of $X$ is called $\mathfrak{c}$-open if its complement is soft $\mathfrak{c}$-closed. By Lemma 4.2 (1), (4) and (5), one can see that the family of all $\mathfrak{c}$-open sets over $X$ forms a soft topology on $X$ and it is called a cluster soft topology or shortly a soft $\mathfrak{c}$-topology. We denote a cluster soft topology on $X$ by $\tilde{\mathcal{T}}_{\mathfrak{c}}(\tilde{I})$ or simply $\tilde{\mathcal{T}}_{\mathfrak{c}}$ when no confusion caused. Each soft open set is soft $\mathfrak{c}$-open but not the reverse (by Lemma 4.2 (3)). This proves that $(X, \tilde{\mathcal{T}}_{\mathfrak{c}}, \Omega)$ is finer than $(X, \tilde{\mathcal{T}}, \Omega)$.*

Notice that crisp (general) cluster topologies are called ideal topologies (see Jankovic and Hamlett (1990)).

The Lemma 11 and Theorem 6 in Azzam et al. (2022) guarantee that the cluster soft topology is equivalent soft $\tilde{\mathcal{T}}^*$-topology constructed differently in Kandil et al. (2014).

**Theorem 4.5** *Let $(X, \tilde{\mathcal{T}}, \Omega)$ be a soft topological space and $\tilde{I} \in I_\Omega(X)$. Then*

$$\mathcal{B} = \{(R, \Omega) - (A, \Omega): (R, \Omega) \in \tilde{\mathcal{T}}, (A, \Omega) \in \tilde{I}\}$$

forms a soft base for the cluster soft topology $\tilde{\mathcal{T}}_{\mathfrak{c}}$.

***Proof***. First we show that $\mathcal{B}$ covers $\tilde{X}$. Let $x_\alpha \in P_\Omega(X)$ and let $(U,\Omega) \in \tilde{\mathcal{T}}_{\mathfrak{c}}$ that contains $x_\alpha$. Then $\tilde{X} - (U,\Omega)$ is soft $\mathfrak{c}$-closed and so $\mathfrak{c}(\tilde{X} - (U,\Omega)) \tilde{\subseteq} \tilde{X} - (U,\Omega)$. This implies that $(U,\Omega) \tilde{\subseteq} \tilde{X} - \mathfrak{c}(\tilde{X} - (U,\Omega))$ and $x_\alpha \notin \mathfrak{c}(\tilde{X} - (U,\Omega))$. Then there exists $(W,\Omega) \in \tilde{\mathcal{T}}(x_\alpha)$ such that $(W,\Omega) \tilde{\cap} [\tilde{X} - (U,\Omega)] \in \tilde{I}$. If we set $(A,\Omega) = (W,\Omega) \tilde{\cap} [\tilde{X} - (U,\Omega)]$, then we will have $x_\alpha \in (W,\Omega) - (A,\Omega) \tilde{\subseteq} (U,\Omega)$.

The proof of Theorem 3.4 in Kandil et al. (2014) showed that $\mathcal{B}$ is closed under finite intersections. Thus, the result is proved.

**Lemma 4.6** *If* $(X,\tilde{\mathcal{T}}_{\mathfrak{c}}(\tilde{I}),\Omega)$ *is a cluster soft topological space, then*

$$\tilde{\mathcal{T}}_{\mathfrak{c}}(\tilde{I})(\alpha) = \tilde{\mathcal{T}}(\alpha)(\tilde{I}(\alpha)),$$

*where* $\tilde{\mathcal{T}}(\alpha)(\tilde{I}(\alpha))$ *is a crisp cluster (ideal) topology with respect to the crisp topoloft* $\tilde{\mathcal{T}}(\alpha)$ *and crisp ideal* $\tilde{I}(\alpha) = \{A(\alpha): (A,\Omega) \in \tilde{I}\}$*, and* $\tilde{\mathcal{T}}_{\mathfrak{c}}(\tilde{I})(\alpha) = \{R(\alpha): (R,\Omega) \in \tilde{\mathcal{T}}_{\mathfrak{c}}(\tilde{I})\}$ *for each* $\alpha \in \Omega$.

***Proof***. Let $\alpha \in \Omega$ and let $R(\alpha) \in \tilde{\mathcal{T}}_{\mathfrak{c}}(\tilde{I})(\alpha)$. Then $(R,\Omega) \in \tilde{\mathcal{T}}_{\mathfrak{c}}(\tilde{I})$. By Theorem 4.5,

$$(R,\Omega) = \tilde{\bigcup}_{j\in J}\left[(G_j,\Omega) - (A_j,\Omega)\right],$$

where $(G_j,\Omega) \in \tilde{\mathcal{T}}$ and $(A_j,\Omega) \in \tilde{I}$. This means that

$$R(\alpha) = \tilde{\bigcup}_{j\in J}[G_j(\alpha) - A_j(\alpha)],$$

where $G_j(\alpha) \in \tilde{\mathcal{T}}(\alpha)$ and $A_j(\alpha) \in \tilde{I}(\alpha)$. Since $\tilde{\mathcal{T}}(\alpha)$ and $\tilde{I}(\alpha)$ are crips topology and ideal for each $\alpha$ (by Lemmas 2.5 and 2.10), therefore, by Theorem 3.1 in Jankovic and Hamlett (1990), each $G_j(\alpha) - A_j(\alpha)$ is basic open set in crisp cluster topology $\tilde{\mathcal{T}}(\alpha)(\tilde{I}(\alpha))$. Thus, $R(\alpha)$ is open in $\tilde{\mathcal{T}}(\alpha)(\tilde{I}(\alpha))$ and so, $\tilde{\mathcal{T}}_{\mathfrak{c}}(\tilde{I})(\alpha) \subseteq \tilde{\mathcal{T}}(\alpha)(\tilde{I}(\alpha))$. The reverse of the inclusion can be proved by a similar technique. Hence, $\tilde{\mathcal{T}}_{\mathfrak{c}}(\tilde{I})(\alpha) = \tilde{\mathcal{T}}(\alpha)(\tilde{I}(\alpha))$.

Nevertheless, there might exist a family of soft sets (need not be a soft topology) and a soft ideal on a set $X$ for which their crisp parts generate cluster topologies, as shown in the following illustration:

**Example 4.7** *Let* $X = \{x,y,z\}$ *and* $\Omega = \{\alpha,\beta\}$*. Consider the soft ideal* $\tilde{I} = \{\tilde{\Phi}, (A_1,\Omega), (A_2,\Omega), \ldots, (A_{15},\Omega)\}$*, where*

$(A_1,\Omega) = \{(\alpha,\emptyset),(\beta,\{x\})\}$,

$(A_2,\Omega) = \{(\alpha,\emptyset),(\beta,\{z\})\}$,

$(A_3,\Omega) = \{(\alpha,\emptyset),(\beta,\{x,z\})\}$,

$(A_4,\Omega) = \{(\alpha,\{y\}),(\beta,\emptyset)\}$,

$(A_5,\Omega) = \{(\alpha,\{y\}),(\beta,\{x\})\}$,

$(A_6,\Omega) = \{(\alpha,\{y\}),(\beta,\{z\})\}$,

$(A_7,\Omega) = \{(\alpha,\{y\}),(\beta,\{x,z\})\}$,

$(A_8,\Omega) = \{(\alpha,\{z\}),(\beta,\emptyset)\}$,

$(A_9, \Omega) = \{(\alpha, \{z\}), (\beta, \{x\})\}$,

$(A_{10}, \Omega) = \{(\alpha, \{z\}), (\beta, \{z\})\}$,

$(A_{11}, \Omega) = \{(\alpha, \{z\}), (\beta, \{x, z\})\}$,

$(A_{12}, \Omega) = \{(\alpha, \{y, z\}), (\beta, \emptyset)\}$,

$(A_{13}, \Omega) = \{(\alpha, \{y, z\}), (\beta, \{x\})\}$,

$(A_{14}, \Omega) = \{(\alpha, \{y, z\}), (\beta, \{z\})\}, and$

$(A_{15}, \Omega) = \{(\alpha, \{y, z\}), (\beta, \{x, z\})\}$.

Let $T = \{\widetilde{\Phi}, (R_1, \Omega), (R_2, \Omega), (R_3, \Omega), (R_4, \Omega), \tilde{X}\}$ be family of soft sets, where

$(R_1, \Omega) = \{(\alpha, \{x\}), (\beta, \{y\})\}$,

$(R_2, \Omega) = \{(\alpha, \{x, y\}), (\beta, \{x, y\})\}$,

$(R_3, \Omega) = \{(\alpha, \{x, y\}), (\beta, \{y, z\})\}, and$

$(R_4, \Omega) = \{(\alpha, \{x, z\}), (\beta, \{x, y\})\}$.

One can easily conclude that $T$ is not a soft topology, which means that we cannot generate a cluster soft topology from $T$ and $\tilde{I}$. On the other hand, both of $T(\alpha) = \{\emptyset, \{x\}, \{x, y\}, \{x, z\}, X\}$ and $T(\beta) = \{\emptyset, \{y\}, \{x, y\}, \{y, z\}, X\}$ define crisp topologies with respect to the respective ideals $\tilde{I}(\alpha) = \{\emptyset, \{y\}, \{z\}, \{y, z\}\}$ and $\tilde{I}(\beta) = \{\emptyset, \{x\}, \{z\}, \{x, z\}\}$ such that $T(\alpha) = T_{\mathfrak{c}}(\alpha)$ and $T(\beta) = T_{\mathfrak{c}}(\beta)$.

However, one can always construct a cluster soft topology by the following scheme:

**Remark 4.8** *Let $(X, \tilde{\mathcal{T}}, \Omega)$ be a soft topological space and let $\tilde{I}$ be a soft ideal. By Lemma 2.9, $(\Omega \times X, \mathcal{T})$ and $I = \{Gr(A): (A, \Omega) \in \tilde{I}\}$ are respectively a crisp topological space and an ideal on $\Omega \times X$. Furthermore, $(X, \tilde{\mathcal{T}}_{\mathfrak{c}}, \Omega)$ forms a cluster soft topological space in which $\mathfrak{c}_{(\tilde{\mathcal{T}}, \tilde{I})}(R, \Omega)$ is a soft set given by a set-valued mapping which graph is equal to $(D_{(\mathcal{T}, I)}(Gr(R))$, where $D_{(\mathcal{T}, I)}(A) = \{(\alpha, x) \in \Omega \times X: A \cap U \notin I, (\alpha, x) \in U \in \mathcal{T}\}, A \subset \Omega \times X$* (see Definition 2.2, Jankovic and Hamlett (1990)). *This means that many results concerning cluster soft topologies can be derived from corresponding results from the theory of ideal topological spaces.*

Here, we shall acknowledge that the one-to-one correspondence between soft and crisp topologies mentioned in the preceding remark was suggested by one of the reviewers.

**Proposition 4.9** *Let $(X, \tilde{\mathcal{T}}, \Omega)$ be a soft topological space and $\tilde{I}, \tilde{J} \in I_{\Omega}(X)$. For any $(R, \Omega) \in S_{\Omega}(X)$, we have*

$$\mathfrak{c}_{(\tilde{\mathcal{T}}, \tilde{I} \tilde{\Cup} \tilde{J})}(R, \Omega) = \mathfrak{c}_{(\tilde{\mathcal{T}}_{\mathfrak{c}}(\tilde{J}), \tilde{I})}(R, \Omega) \widetilde{\bigcap} \mathfrak{c}_{(\tilde{\mathcal{T}}_{\mathfrak{c}}(\tilde{I}), \tilde{J})}(R, \Omega).$$

***Proof*.** Given $x_\alpha \in P_{\Omega}(X)$ and suppose $x_\alpha \notin \mathfrak{c}_{(\tilde{\mathcal{T}}, \tilde{I} \tilde{\Cup} \tilde{J})}(R, \Omega)$. Then there is $(U, \Omega) \in \tilde{\mathcal{T}}(x_\alpha)$ such that $(R, \Omega) \tilde{\cap} (U, \Omega) \in \tilde{I} \tilde{\Cup} \tilde{J}$. Set $(R, \Omega) \tilde{\cap} (U, \Omega) = (A, \Omega) \tilde{\cup} (B, \Omega)$ for some $(A, \Omega), (B, \Omega)$ in $\tilde{I}, \tilde{J}$ respectively. Since $\tilde{I}$ is closed under soft subsets, so we can assume that $(A, \Omega) \tilde{\cap} (B, \Omega) = \widetilde{\Phi}$. Therefore, $[(R, \Omega) \tilde{\cap} (U, \Omega)] - (A, \Omega) = (B, \Omega)$ and $[(R, \Omega) \tilde{\cap} (U, \Omega)] - (B, \Omega) = (A, \Omega)$ . This means $[(R,\Omega) \tilde{\cap} (U,\Omega)] - (A,\Omega) \in \tilde{J}$ and $[(R,\Omega) \tilde{\cap} (U,\Omega)] - (B,\Omega) \in \tilde{I}$. Then we have $x_\alpha \notin \mathfrak{c}_{(\tilde{\mathcal{T}}_{\mathfrak{c}}(\tilde{I}),\tilde{J})}(R,\Omega)$ or $x_\alpha \notin \mathfrak{c}_{(\tilde{\mathcal{T}}_{\mathfrak{c}}(\tilde{J}),\tilde{I})}(R,\Omega)$ as either $x_\alpha \in \tilde{I}$ or $x_\alpha \in \tilde{J}$. Thus,

$$\mathfrak{c}_{(\tilde{\mathcal{T}}_{\mathfrak{c}}(\tilde{J}),\tilde{I})}(R,\Omega) \tilde{\cap} \mathfrak{c}_{(\tilde{\mathcal{T}}_{\mathfrak{c}}(\tilde{I}),\tilde{J})}(R,\Omega) \tilde{\subseteq} \mathfrak{c}_{(\tilde{\mathcal{T}},\tilde{I}\tilde{\Cup}\tilde{J})}(R,\Omega).$$

Conversely, if $x_\alpha \notin \mathfrak{c}_{(\tilde{\mathcal{T}}_{\mathfrak{c}}(\tilde{I}),\tilde{J})}(R,\Omega)$, then there exist $(V,\Omega) \in \tilde{\mathcal{T}}(x_\alpha)$ and $(A,\Omega) \in \tilde{I}$ such that $[(V,\Omega) - (A,\Omega)] \tilde{\cap} (R,\Omega) \in \tilde{J}$. Since $\tilde{I}$ is closed under soft subsets, we let $(A,\Omega) \tilde{\subseteq} (R,\Omega)$. If $(B,\Omega) = [(V,\Omega) - (A,\Omega)] \tilde{\cap} (R,\Omega)$, then $(V,\Omega) \tilde{\cap} (R,\Omega) = (A,\Omega) \tilde{\cup} (B,\Omega) \in \tilde{I} \tilde{\Cup} \tilde{J}$ and hence $x_\alpha \notin \mathfrak{c}_{(\tilde{\mathcal{T}},\tilde{I}\tilde{\Cup}\tilde{J})}(R,\Omega)$. This shows that $\mathfrak{c}_{(\tilde{\mathcal{T}},\tilde{I}\tilde{\Cup}\tilde{J})}(R,\Omega) \tilde{\subseteq} \mathfrak{c}_{(\tilde{\mathcal{T}}_{\mathfrak{c}}(\tilde{I}),\tilde{J})}(R,\Omega)$. Symmetrically, we can obtain that $\mathfrak{c}_{(\tilde{\mathcal{T}},\tilde{I}\tilde{\Cup}\tilde{J})}(R,\Omega) \tilde{\subseteq} \mathfrak{c}_{(\tilde{\mathcal{T}}_{\mathfrak{c}}(\tilde{J}),\tilde{I})}(R,\Omega)$. In conclusion, we get $\mathfrak{c}_{(\tilde{\mathcal{T}},\tilde{I}\tilde{\Cup}\tilde{J})}(R,\Omega) \tilde{\subseteq} \mathfrak{c}_{(\tilde{\mathcal{T}}_{\mathfrak{c}}(\tilde{I}),\tilde{J})}(R,\Omega) \tilde{\cap} \mathfrak{c}_{(\tilde{\mathcal{T}}_{\mathfrak{c}}(\tilde{J}),\tilde{I})}(R,\Omega)$. The proof ends.

The next result illustrates that by constructing the cluster soft topology twice, you will get the first obtained cluster soft topology.

**Theorem 4.10** *Let $(X,\tilde{\mathcal{T}},\Omega)$ be a soft topological space and $\tilde{I} \in I_\Omega(X)$. Then $\tilde{\mathcal{T}}_{\mathfrak{cc}} = \tilde{\mathcal{T}}_{\mathfrak{c}}$.*

***Proof.*** By assuming $\tilde{I} = \tilde{J}$ in Proposition 4.9, we obtain that a soft set $(R,\Omega)$ is soft $\mathfrak{c}$-closed iff it is soft $\mathfrak{cc}$-closed. This implies that $\tilde{\mathcal{T}}_{\mathfrak{cc}} = \tilde{\mathcal{T}}_{\mathfrak{c}}$, where $\tilde{\mathcal{T}}_{\mathfrak{cc}}$ is the cluster soft topology of $\tilde{\mathcal{T}}_{\mathfrak{c}}(\tilde{I})$.

Next, we present a few illustrations highlighting the significance of cluster soft topologies.

**Example 4.11** *Given any soft topological space $(X,\tilde{\mathcal{T}},\Omega)$ and any soft ideal $\tilde{I} \in I_\Omega(X)$. If $\tilde{I}$ is trivial (i.e., $\tilde{I} = \{\tilde{\Phi}\}$), by Remark 3.2 (1), $\mathfrak{c}(R,\Omega) = cl(R,\Omega)$ which implies that $\tilde{\mathcal{T}}_{\mathfrak{c}} = \tilde{\mathcal{T}}$.*

**Example 4.12** *Let $(X,\tilde{\mathcal{T}},\Omega)$ be a soft topological space and let $\tilde{I} \in I_\Omega(X)$. If $\tilde{I} = S_\Omega(X)$, then $\mathfrak{c}(R,\Omega) = \tilde{\Phi}$ which implies that each $(R,\Omega) \in S_\Omega(X)$ is soft $\mathfrak{c}$-closed. Therefore, $\tilde{\mathcal{T}}_{\mathfrak{c}} = \tilde{\mathcal{T}}_{dis}$, where $\tilde{\mathcal{T}}_{dis}$ is the soft discrete topology.*

**Example 4.13** *Let $(X,\tilde{\mathcal{T}}_{ind},\Omega)$ be the indiscrete soft topological space and let $\tilde{I} \in I_\Omega(X)$, where $\tilde{I} = \{(A,\Omega): (A,\Omega) \in S_\Omega(X), (A,\Omega) \text{ is finite}\}$. If $(R,\Omega) \in \tilde{I}$, then $\mathfrak{c}(R,\Omega) = \tilde{\Phi}$. If $(R,\Omega) \notin \tilde{I}$, then $\mathfrak{c}(R,\Omega) = \tilde{X}$. Consequently, each finite soft set is soft $\mathfrak{c}$-closed together with $\tilde{X}$. Therefore, $\tilde{\mathcal{T}}_{\mathfrak{c}} = \tilde{\mathcal{T}}_{cof}$, where $\tilde{\mathcal{T}}_{cof}$ is the soft co-finite topology (c.f., Theorem 5.1).*

**Example 4.14** *Let $(X,\tilde{\mathcal{T}}_{ind},\Omega)$ be the indiscrete soft topological space, $x_\alpha \in S_\Omega(X)$, and $\tilde{I} \in I_\Omega(X)$. Suppose $\tilde{I} = \{(A,\Omega): (A,\Omega) \in S_{\Omega(X)}, x_\alpha \notin (A,\Omega)\}$. If $(R,\Omega) \in \tilde{I}$, then $\mathfrak{c}(R,\Omega) = \tilde{\Phi}$. If $(R,\Omega) \notin \tilde{I}$, then $\mathfrak{c}(R,\Omega) = \tilde{X}$. Therefore, each soft set excluding $x_\alpha$ is soft $\mathfrak{c}$-closed together with $\tilde{X}$. Therefore, $\tilde{\mathcal{T}}_{\mathfrak{c}} = \tilde{\mathcal{T}}_{inc}$, where $\tilde{\mathcal{T}}_{inc} = \{(A,\Omega): (A,\Omega) \in S_\Omega(X), x_\alpha \in (A,\Omega)\} \tilde{\cup} \{\tilde{\Phi}\}$, (it is called included soft point topology in Example 2 in Al Ghour and Ameen (2022)).*

# 5 Characterizations of soft $\mathfrak{c}$-closed sets

In this section, we characterize soft $\mathfrak{c}$-closed sets in terms of some other classes of soft sets when the underlying soft topology or the related soft ideal possesses certain properties.

**Theorem 5.1** *Let $\tilde{\mathcal{T}} \in T_{\Omega}(X)$ and $\tilde{I} = \{(A,\Omega): (A,\Omega) \in S_{\Omega}(X), (A,\Omega)$ is finite$\}$. Then each soft $\mathfrak{c}$-closed set is soft closed iff each finite soft set is soft closed.*

***Proof***. The first part is easy as each finite soft set in $\tilde{I}$, by Lemma 4.2 (2), all finite soft sets are $\mathfrak{c}$-closed and so they are soft closed by the assumption.

Conversely, suppose each finite soft set is soft closed and $x_\alpha \in P_\Omega(X)$. Let $(R,\Omega)$ be a soft $\mathfrak{c}$-closed set. If $x_\alpha \notin (R,\Omega)$, then by Remark 3.2 (2), $x_\alpha \notin \mathcal{D}(R,\Omega)$. Therefore, there exists $(U,\Omega) \in \tilde{\mathcal{T}}(x_\alpha)$ such that $(R,\Omega) \tilde{\cap} (U,\Omega)$ is a finite soft set. Set $(Q,\Omega) = (R,\Omega) \tilde{\cap} (U,\Omega)$. Since $x_\alpha \notin (Q,\Omega)$, so the soft set $(V,\Omega) = (U,\Omega) - (Q,\Omega) \in \tilde{\mathcal{T}}(x_\alpha)$ and $(R,\Omega) \tilde{\cap} (V,\Omega) = \tilde{\Phi}$. This means that $x_\alpha \notin cl(R,\Omega)$. Hence, $(R,\Omega)$ is a soft closed set.

**Corollary 5.2** *Let $\tilde{\mathcal{T}} \in T_{\Omega}(X)$ and $\tilde{I} = \{(A,\Omega): (A,\Omega) \in S_{\Omega}(X), (A,\Omega)$ is finite$\}$. Then $\tilde{\mathcal{T}} = \tilde{\mathcal{T}}_{\mathfrak{c}}$ iff the complement of each finite soft set is soft open.*

**Definition 5.3** *Given $\tilde{\mathcal{T}} \in T_{\Omega}(X)$ and $\tilde{I} \in I_{\Omega}(X)$. Then $\tilde{I}$ is called a soft adherent ideal if $(R,\Omega) - \mathfrak{c}(R,\Omega) \in \tilde{I}$ for each $(R,\Omega) \in S_{\Omega}(X)$.*

The next example illustrates that a soft ideal of finite soft sets need not be soft adherent.

**Example 5.4** *Let $(X,\tilde{\mathcal{T}}_{dis},\Omega)$ be the discrete soft topological space, where $X$ is infinite, and let $\tilde{I}$ be a soft ideal defined by $\tilde{I} = \{(A,\Omega): (A,\Omega) \in S_{\Omega}(X), (A,\Omega)$ is finite$\}$. Then $\mathfrak{c}(R,\Omega) = \tilde{\Phi}$ for all $(R,\Omega) \in S_{\Omega}(X)$. Therefore, for each infinite $(R,\Omega) \in S_{\Omega}(X)$, $(A,\Omega) - \mathfrak{c}(A,\Omega) \notin \tilde{I}$. Thus, $\tilde{I}$ is not soft adherent.*

However, the following result informs us that a soft $\sigma$-ideal in certain soft topological spaces is soft adherent.

**Theorem 5.5** *If $\tilde{\mathcal{T}} \in T_{\Omega}(X)$ has a countable soft base, then each soft $\sigma$-ideal on $X$ is soft adherent.*

***Proof***. Let $\tilde{I}$ be a soft $\sigma$-ideal on $X$ and let $(R,\Omega) \in S_{\Omega}(X)$. For each $x_\alpha \in (R,\Omega) - \mathfrak{c}(R,\Omega)$ we associate a $(U_{x_\alpha},\Omega) \in \tilde{\mathcal{T}}(x_\alpha)$ such that $(R,\Omega) \tilde{\cap} (U_{x_\alpha},\Omega) \in \tilde{I}$. Set $\mathcal{U} = \{(U_{x_\alpha},\Omega): x_\alpha \in (R,\Omega) - \mathfrak{c}(R,\Omega)\}$. Since $\tilde{\mathcal{T}}$ has a countable soft base, so there is a countable soft subset $(D,\Omega)$ of $(R,\Omega) - \mathfrak{c}(R,\Omega)$ such that $\tilde{\cup}\mathcal{U} = \tilde{\cup}_{x_\alpha \in (D,\Omega)} (U_{x_\alpha},\Omega)$. Therefore, we have

$$(R,\Omega) - \mathfrak{c}(R,\Omega) \tilde{\sqsubseteq} \tilde{\bigcup}_{x_\alpha \in (R,\Omega)-\mathfrak{c}(R,\Omega)} [(R,\Omega) \tilde{\cap} (U_{x_\alpha},\Omega)] = \tilde{\bigcup}_{x_\alpha \in (D,\Omega)} [(R,\Omega) \tilde{\cap} (U_{x_\alpha},\Omega)].$$

Since $\tilde{\bigcup}_{x_\alpha \in (D,\Omega)} [(R,\Omega) \tilde{\cap} (U_{x_\alpha},\Omega)] \in \tilde{I}$, then $(R,\Omega) - \mathfrak{c}(R,\Omega) \in \tilde{I}$ and hence $\tilde{I}$ is soft adherent.

**Theorem 5.6** *Let $\tilde{\mathcal{T}} \in T_{\Omega}(X)$ and let $\tilde{I} \in I_{\Omega}(X)$ be soft adherent. A soft set $(R,\Omega)$ over $X$ is soft $\mathfrak{c}$-closed iff it can be written as a disjoint union of a soft closed set and an element in $\tilde{I}$.*

***Proof***. Suppose $(R,\Omega)$ is soft $\mathfrak{c}$-closed. Set $(Q,\Omega) = \mathfrak{c}(R,\Omega)$. By Lemma 3.5 (2), $(Q,\Omega)$ is soft closed. If $(A,\Omega) = (Q,\Omega) - \mathfrak{c}(R,\Omega)$, by assumption, $(A,\Omega) \in \tilde{I}$. Thus, $(R,\Omega) = (Q,\Omega) \tilde{\cup} (A,\Omega)$ is the required form.

Conversely, suppose such that representation exists. Since each element $(A,\Omega)$ in $\tilde{I}$ is soft $\mathfrak{c}$-closed and each soft closed set $(Q,\Omega)$ is soft $\mathfrak{c}$-closed, so $(Q,\Omega) \tilde{\cup} (A,\Omega)$ is soft $\mathfrak{c}$-closed (*c.f.*, Theorem 4.3 in Kandil et al. (2014)).

The following is a direct consequences of the above result:

**Corollary 5.7** *Let $\tilde{\mathcal{T}} \in T_\Omega(X)$ and let $\tilde{I} \in I_\Omega(X)$ be soft adherent. A soft set $(W,\Omega)$ is soft $\mathfrak{c}$-open iff it has the form $(W,\Omega) = (U,\Omega) - (A,\Omega)$, where $(U,\Omega) \in \tilde{\mathcal{T}}$ and $(A,\Omega) \in \tilde{I}$.*

**Theorem 5.8** *Let $(X,\tilde{\mathcal{T}},\Omega)$ be a soft topological space of a countable soft base $\mathcal{B}$ and let $\tilde{I}$ be a soft $\sigma$-ideal on $\tilde{X}$. Then $\tilde{\mathcal{T}}_\mathfrak{c} = \{(R,\Omega) - (A,\Omega): (R,\Omega) \in \tilde{\mathcal{T}}, (A,\Omega) \in \tilde{I}\}$.*

***Proof***. Let $\mathcal{B} = \{(R,\Omega) - (A,\Omega): (R,\Omega) \in \tilde{\mathcal{T}}, (A,\Omega) \in \tilde{I}\}$. We first show that $\tilde{\mathcal{T}}_\mathfrak{c} \tilde{\subseteq} \mathcal{B}$. Let $(R,\Omega) \in \tilde{\mathcal{T}}_\mathfrak{c}$. Since $(X,\tilde{\mathcal{T}},\Omega)$ has a countable soft base, by Theorem 5.5, $\tilde{I}$ is a soft adherent ideal. Then Corollary 5.7 guarantees that $(R,\Omega) = (U,\Omega) - (A,\Omega)$ for some $(U,\Omega) \in \tilde{\mathcal{T}}$ and $(A,\Omega) \in \tilde{I}$. This implies by Theorem 4.5 that $(R,\Omega) \in \mathcal{B}$. On the other hand, if $(R,\Omega) \in \mathcal{B}$, then $(R,\Omega) = (U,\Omega) - (A,\Omega)$, where $(U,\Omega) \in \tilde{\mathcal{T}}$ and $(A,\Omega) \in \tilde{I}$. By the same reason above, $(R,\Omega) \in \tilde{\mathcal{T}}_\mathfrak{c}$ and thus, $\mathcal{B} \tilde{\subseteq} \tilde{\mathcal{T}}_\mathfrak{c}$. Hence, $\tilde{\mathcal{T}}_\mathfrak{c} = \mathcal{B}$.

**Definition 5.9** *Let $(R,\Omega) \in S_\Omega(X)$, $\tilde{\mathcal{T}} \in T_\Omega(X)$, and $\tilde{I} \in I_\Omega(X)$. Then $(R,\Omega)$ is said to be*

1. soft $\mathfrak{c}$-crowded if $(R,\Omega) \tilde{\subseteq} \mathfrak{c}(R,\Omega)$.
2. soft $\mathfrak{c}$-regular if $\mathfrak{c}(R,\Omega) = (R,\Omega)$.

**Theorem 5.10** *Let $\tilde{\mathcal{T}} \in T_\Omega(X)$ and $\tilde{I} \in I_\Omega(X)$. If $\tilde{I}$ is a soft adherent-ideal, then each soft set can be represented as a disjoint union of a soft $\mathfrak{c}$-crowded set and an element in $\tilde{I}$.*

***Proof***. Given a soft set $(R,\Omega)$. Consider the following decomposition

$$(R,\Omega) = [(R,\Omega) \tilde{\cap} \mathfrak{c}(R,\Omega)] \tilde{\cup} [(R,\Omega) - \mathfrak{c}(R,\Omega)]. \quad (1)$$

Since $\tilde{I}$ is soft adherent, $(R,\Omega) - \mathfrak{c}(R,\Omega) \in \tilde{I}$. It remains to show that $(R,\Omega) \tilde{\cap} \mathfrak{c}(R,\Omega)$ is a soft $\mathfrak{c}$-crowded set. Now, by Proposition 3.3 (1) and (4), we have

$$(R,\Omega) \tilde{\cap} \mathfrak{c}(R,\Omega) \tilde{\subseteq} \mathfrak{c}(R,\Omega) = \mathfrak{c}[(R,\Omega) \tilde{\cap} \mathfrak{c}(R,\Omega)] \tilde{\cup} \mathfrak{c}[(R,\Omega) - \mathfrak{c}(R,\Omega)] = \mathfrak{c}[(R,\Omega) \tilde{\cap} \mathfrak{c}(R,\Omega)].$$

Therefore, $(R,\Omega) \tilde{\cap} \mathfrak{c}(R,\Omega)$ is a soft $\mathfrak{c}$-crowded set. Hence the proof.

**Lemma 5.11** *Let $(R,\Omega),(S,\Omega) \in S_\Omega(X)$, $\tilde{\mathcal{T}} \in T_\Omega(X)$, and $\tilde{I} \in I_\Omega(X)$. If $(R,\Omega)$ is soft $\mathfrak{c}$-regular and $(S,\Omega)$ is soft $\mathfrak{c}$-closed in $(R,\Omega)$, then $(R,\Omega) - (S,\Omega)$ is soft $\mathfrak{c}$-crowded.*

***Proof***. Given $x_\alpha \in P_\Omega(X)$ and set $(Q,\Omega) = (R,\Omega) - (S,\Omega)$. Suppose $x_\alpha \in (Q,\Omega)$. Since $(Q,\Omega) \tilde{\subseteq} (R,\Omega) \tilde{\subseteq} \mathfrak{c}(R,\Omega)$, then for each $(U,\Omega) \in \tilde{\mathcal{T}}(x_\alpha)$ such that $(U,\Omega) \tilde{\cap} (R,\Omega) \notin \tilde{I}$. Since $x_\alpha \notin (S,\Omega)$ and $(S,\Omega)$ is soft $\mathfrak{c}$-closed, so $x_\alpha \notin \mathfrak{c}(S,\Omega)$. Therefore, there exists $(V,\Omega) \in \tilde{\mathcal{T}}(x_\alpha)$ such that $(V,\Omega) \tilde{\cap} (S,\Omega) \in \tilde{I}$. Set $(W,\Omega) = (U,\Omega) \tilde{\cap} (V,\Omega)$. Also $(W,\Omega) \in \tilde{\mathcal{T}}(x_\alpha)$. Now, we have

$$(R,\Omega) \tilde{\cap} (W,\Omega) = [(Q,\Omega) \tilde{\cap} (W,\Omega)] \tilde{\cup} [(S,\Omega) \tilde{\cap} (W,\Omega)].$$

Since $(W,\Omega) \tilde{\cap} (R,\Omega) \notin \tilde{I}$, we must have $(W,\Omega) \tilde{\cap} (Q,\Omega) \notin \tilde{I}$. Hence, $x_\alpha \in \mathfrak{c}(Q,\Omega)$ and thus $(R,\Omega) - (S,\Omega)$ is soft $\mathfrak{c}$-crowded.

**Theorem 5.12** *Let $\tilde{\mathcal{T}} \in T_\Omega(X)$ and $\tilde{I} \in I_\Omega(X)$. If $\tilde{I}$ is a soft adherent-ideal, then each soft $\mathfrak{c}$-closed set can be uniquely represented as a disjoint union of a soft $\mathfrak{c}$-regular set and an element in $\tilde{I}$.*

***Proof***. Let $(R,\Omega)$ be soft $\mathfrak{c}$-closed. We first prove that such a representation exists and then show that it is unique. Consider,

$$(R,\Omega) = [(R,\Omega) \tilde{\cap} \mathfrak{c}(R,\Omega)] \tilde{\cup} [(R,\Omega) - \mathfrak{c}(R,\Omega)].$$

Since $\tilde{I}$ is soft adherent, $(R,\Omega) - \mathfrak{c}(R,\Omega) \in \tilde{I}$. We now show that $(R,\Omega) \tilde{\cap} \mathfrak{c}(R,\Omega)$ is a soft $\mathfrak{c}$-regular set. Since $(R,\Omega)$ is soft $\mathfrak{c}$-closed, then $\mathfrak{c}(R,\Omega) \tilde{\subseteq} (R,\Omega)$. Therefore, by Proposition 3.3 (2), we have

$$\mathfrak{c}((R,\Omega) \tilde{\cap} \mathfrak{c}(R,\Omega)) \tilde{\subseteq} \mathfrak{c}(R,\Omega) = (R,\Omega) \tilde{\cap} \mathfrak{c}(R,\Omega).$$

The opposite of the inclusion is obtained from Theorem 5.10. Thus, $(R,\Omega) \tilde{\cap} \mathfrak{c}(R,\Omega)$ is a soft $\mathfrak{c}$-regular set.

Suppose $(R,\Omega)$ has two different representations. Namely, $(R,\Omega) = (L,\Omega) \tilde{\cup} (A,\Omega) = (M,\Omega) \tilde{\cup} (B,\Omega)$, where $(L,\Omega),(M,\Omega)$ are disjoint soft $\mathfrak{c}$-regular and $(A,\Omega),(B,\Omega) \in \tilde{I}$ with $(A,\Omega) \tilde{\cap} (B,\Omega) = \tilde{\Phi}$. Since $(L,\Omega) - [(L,\Omega) \tilde{\cap} (M,\Omega)] \tilde{\subseteq} (B,\Omega)$ and $(B,\Omega) \in \tilde{I}$, by Lemma 5.11, we conclude

$$(L,\Omega) - [(L,\Omega) \tilde{\cap} (M,\Omega)] \tilde{\subseteq} \mathfrak{c}[(L,\Omega) - ((L,\Omega) \tilde{\cap} (M,\Omega))] \tilde{\subseteq} \mathfrak{c}(B,\Omega) = \tilde{\Phi}.$$

Thus, $(L,\Omega) = (L,\Omega) \tilde{\cap} (M,\Omega)$. By the same way, one can obtain $(M,\Omega) = (L,\Omega) \tilde{\cap} (M,\Omega)$. Hence, $(L,\Omega) = (M,\Omega)$ and $(A,\Omega) = (B,\Omega)$

**Theorem 5.13** *Let $\tilde{\mathcal{T}} \in T_\Omega(X)$ and let $\tilde{I} \in I_\Omega(X)$ be soft adherent. A soft set $(R,\Omega)$ over $X$ is soft $\mathfrak{c}$-closed iff it can be written as a disjoint union of a soft $\mathfrak{c}$-regular set and an element in $\tilde{I}$.*

***Proof***. The first direction follows from Theorem 5.12. For the converse, if $(R,\Omega) = (Q,\Omega) \tilde{\cup} (A,\Omega)$, where $(Q,\Omega)$ is soft $\mathfrak{c}$-regular and $(A,\Omega) \in \tilde{I}$. Then

$$\mathfrak{c}(R,\Omega) = \mathfrak{c}(Q,\Omega) \tilde{\cup} \mathfrak{c}(A,\Omega) = \mathfrak{c}(Q,\Omega) = (Q,\Omega) \tilde{\subseteq} (Q,\Omega) \tilde{\cup} (A,\Omega) = (R,\Omega).$$

This proves that $(R,\Omega)$ is soft $\mathfrak{c}$-closed.

# Conclusion

In this paper, we have considered a cluster soft point as an extension or a unification of a soft closure and a soft limit point. The class of all cluster soft points of a soft set $(F,\Omega)$ is called a cluster soft set or soft local function of $(F,\Omega)$ in Kandil et al. (2014). The cluster set of crisp points was first studied by Vaidyanathaswamy (1944) and developed by Jankovic and Hamlett (1990). The concept of soft ideal played a significant role in determining the cluster soft set. Then we have defined cluster soft closed sets and shown that the system of all cluster soft open sets, which are the complements of cluster soft closed sets, constitutes a soft topology and is called a cluster soft topology. We have demonstrated that the cluster soft topology $\tilde{\mathcal{T}}_{\mathfrak{c}}$ of the soft topology $\tilde{\mathcal{T}}$ is finer than $\tilde{\mathcal{T}}$ (i.e., $\tilde{\mathcal{T}} \tilde{\subseteq} \tilde{\mathcal{T}}_{\mathfrak{c}}$). In addition, we have established that the cluster soft topology $\tilde{\mathcal{T}}_{\mathfrak{cc}}$ of $\tilde{\mathcal{T}}_{\mathfrak{c}}$ is similar to $\tilde{\mathcal{T}}_{\mathfrak{c}}$. Several examples demonstrated that the cluster soft topologies are among the most natural elements in the lattice of soft topologies on a universal set. We have characterized the cluster soft closed sets in terms of some other types of soft sets. At the end, we shall remark that our results have built on the soft point theory given in (Nazmul and Samanta (2013); Xie (2015)). The obtained results are natural generalizations of those found in (Jankovic and Hamlett (1990); Vaidyanathaswamy (1944)), according to Theorem 1 in Terepeta (2019).

As a piece of future work, one can study the concepts of separation axioms, compactness, connectedness, etc, in cluster soft topological spaces. It is also possible to define and investigate some weak or strong classes of cluster soft open sets after introducing the soft $c$-interior and the soft $c$-closure of a soft set.

## Compliance with ethical standards

**Funding** Not applicable.

**Availability of data andmaterials** Not applicable.

**Conflict of interest** The authors declare that they have no competing interests.

**Human and animal rights statement** This article does not contain any studies with human participants performed by any of the authors.

# References

Acar, U., Koyuncu, F., and Tanay, B. (2010). Soft sets and soft rings. Computers & Mathematics with Applications, 59(11):3458–3463.

Akta¸s, H. and Çağman, N. (2007). Soft sets and soft groups. Information sciences, 177(13):2726– 2735.

Al Ghour, S. and Ameen, Z. A. (2022). On soft submaximal spaces. Heliyon, 8(9):e10574.

Al Ghour, S. and Ameen, Z. A. (2022a). Maximal soft compact and maximal soft connected topologies. Applied Computational Intelligence and Soft Computing, 2022:Article ID 9860015.

Al-shami, T. M. (2022). Soft somewhat open sets: soft separation axioms and medical application to nutrition. Computational and Applied Mathematics, 41(5):1–22.

Al-shami, T. M. and El-Shafei, M. E. (2020). Partial belong relation on soft separation axioms and decision-making problem, two birds with one stone. Soft Computing, 24(7):5377–5387.

Al-shami, T. M. and Kocinac, L. D. (2019). The equivalence between the enriched and extended soft topologies. Appl. Comput. Math, 18(2):149–162.

Alcantud, J. C. R. (2020). Soft open bases and a novel construction of soft topologies from bases for topologies. Mathematics, 8(5):672.

Ali, M. I. (2011). A note on soft sets, rough soft sets and fuzzy soft sets. Applied Soft Computing, 11(4):3329–3332.

Ali, M. I., Feng, F., Liu, X., Min, W. K., and Shabir, M. (2009). On some new operations in soft set theory. Computers & Mathematics with Applications, 57(9):1547–1553.

Ameen, Z. A. and Al Ghour, S. (2022). Extensions of soft topologies. Filomat, 36(15):5279–5287

Ameen, Z. A. and Al Ghour, S. (2022a). Minimal soft topologies. New Mathematics and Natural Computation, 19(1), 19-31

Asaad, B. A. (2017). Results on soft extremally disconnectedness of soft topological spaces. J. Math. Computer Sci, 17:448–464.

Aygünoğlu, A. and Aygün, H. (2012). Some notes on soft topological spaces. Neural computing and Applications, 21(1):113–119.

Azzam, A., Ameen, Z. A., Al-shami, T. M., and El-Shafei, M. E. (2022). Generating soft topologies via soft set operators. Symmetry, 14(5):914.

Bayramov, S. and Gunduz, C. (2018). A new approach to separability and compactness in soft topological spaces. TWMS Journal of Pure and Applied Mathematics, 9(21):82–93.

Çağman, N., Karataş, S., and Enginoglu, S. (2011). Soft topology. Computers & Mathematics with Applications, 62(1):351–358.

Dalkılıc¸, O. (2021). A novel approach to soft set theory in decision-making under uncertainty. International Journal of Computer Mathematics, 98(10):1935–1945.

Dalkılıç, O. (2022a). On topological structures of virtual fuzzy parametrized fuzzy soft sets. Complex & Intelligent Systems, 8(1):337–348.

Dalkılıç, O. (2022b). Two novel approaches that reduce the effectiveness of the decision maker in decision making under uncertainty environments. Iranian Journal of Fuzzy Systems, 19(2):105–117.

Dalkılıç, O. and Demirtaş, N. (2022). Algorithms for covid-19 outbreak using soft set theory: estimation and application. Soft Computing, pages 1–9.

Das, S. and Samanta, S. (2013). Soft metric. Ann. Fuzzy Math. Inform, 6(1):77–94.

Dubois, D. and Prade, H. (1990). Rough fuzzy sets and fuzzy rough sets. International Journal of General System, 17(2-3):191–209.

Janković, D. and Hamlett, T. (1990). New topologies from old via ideals. The American mathematical monthly, 97(4):295–310.

Kandil, A., AE Tantawy, O., A El-Sheikh, S., and M Abd El-latif, A. (2014). Soft ideal theory soft local function and generated soft topological spaces. Applied Mathematics & Information Sciences, 8(4):1595–1603.

Lin, F. (2013). Soft connected spaces and soft paracompact spaces. International Journal of Mathematical and Computational Sciences, 7(2):277–283.

Liu, X., Feng, F., Wang, Q., Yager, R. R., Fujita, H., and Alcantud, J. C. R. (2021). Mining temporal association rules with temporal soft sets. Journal of Mathematics, 2021.

Maji, P. K., Biswas, R., and Roy, A. (2001). Fuzzy soft sets. Journal of Fuzzy Mathematics, 9(3):589–602.

Maji, P. K., Biswas, R., and Roy, A. R. (2003). Soft set theory. Computers & Mathematics with Applications, 45(4-5):555–562.

Matejdes, M. (2021). Methodological remarks on soft topology. Soft computing, 25(5):4149– 4156.

Matejdes, M. (2016). Soft topological space and topology on the cartesian product. Hacettepe Journal of Mathematics and Statistics, 45(4):1091–1100.

Molodtsov, D. (1999). Soft set theory-first results. Computers & Mathematics with Applications, 37(4-5):19–31.

Molodtsov, D. (2004). The theory of soft sets. URSS Publishers, Moscow.

Nazmul, S. and Samanta, S. (2013). Neighbourhood properties of soft topological spaces. Ann. Fuzzy Math. Inform, 6(1):1–15.

Pawlak, Z. (1982). Rough sets. International journal of computer & information sciences, 11(5):341–356.

Santos-García, G. and Alcantud, J. C. R. (2023). Ranked soft sets. Expert Systems, page e13231.

Sardar, S. K. and Gupta, S. (2013). Soft category theory-an introduction. Journal of Hyperstructures, 2(2).

Shabir, M. and Naz, M. (2011). On soft topological spaces. Computers & Mathematics with Applications, 61(7):1786–1799.

Terepeta, M. (2019). On separating axioms and similarity of soft topological spaces. Soft Computing, 23(3):1049–1057.

Vaidyanathaswamy, R. (1944). The localisation theory in set-topology. In Proceedings of the Indian Academy of Sciences-Section A, volume 20, pages 51–61. Springer India.

Xie, N. (2015). Soft points and the structure of soft topological spaces. Ann. Fuzzy Math. Inform, 10(2):309–322.

Zadeh, L. (1965). Fuzzy sets. Information and Control, 8(3):338–353.

Zhan, J. and Alcantud, J. C. R. (2019). A survey of parameter reduction of soft sets and corresponding algorithms. Artificial Intelligence Review, 52(3):1839–1872.